\documentclass[11pt]{article}
\usepackage{amsmath,amssymb,amsfonts}
\usepackage[latin1]{inputenc}
\usepackage{epsfig}
\numberwithin{equation}{section}
\newtheorem{thm}{Theorem}[section]
\newtheorem{aadef}[thm]{Definition}
\newtheorem{alem}[thm]{Lemma}
\newtheorem{aprop}[thm]{Proposition}

\newtheorem{arem}[thm]{Remark}
\newenvironment{adem}[1][]%
   {\ \\ {\bf Proof #1~: }}%
   {\hfill\mbox{\rule{2 true mm}{3 true mm}}\vskip 2 ex\noindent}
   {\ \\ {\bf Example #1~: }}%
   {\hfill\mbox{\rule{2 true mm}{3 true mm}}\vskip 2 ex\noindent}

\newcommand{\E}{{\mathbb E}}

\newcommand{\N}{{\mathbb N}}
\newcommand{\R}{{\mathbb R}}

\title{Equivalence of the Poincaré inequality with a transport-chi-square inequality in dimension one}
\author{B.Jourdain\thanks{Universit\'e Paris-Est, CERMICS, 6-8 av Blaise Pascal, Cit\'e
    Descartes, Champs sur Marne, 77455 Marne-la-Vall\'ee Cedex 2, France -
    e-mail : jourdain@cermics.enpc.fr }}

\begin{document}
\maketitle  
\begin{abstract}
  In this paper, we prove that, in dimension one, the Poincaré inequality is equivalent to a new transport-chi-square inequality linking the square of the quadratic Wasserstein distance with the chi-square pseudo-distance. We also check tensorization of this transport-chi-square inequality.
 \end{abstract}
For $q\geq 1$, the Wasserstein distance with index $q$ between two probability measures $\mu$ and $\nu$ on $\R^d$ is denoted by
\begin{equation}
   W_q^q(\mu,\nu)=\inf_{\gamma<^{\mu}_{\nu}}\int_{\R^d\times\R^d}|x-y|^qd\gamma(x,y)\label{wass}
\end{equation}
where the infimum is taken over all probability measures $\gamma$ on $\R^d\times\R^d$ with respective marginals $\mu$ and $\nu$.
We also introduce the relative entropy and the chi-square pseudo distance
\begin{align*}
   H(\nu|\mu)&=\begin{cases}\int_{\R^d} \ln\left(\frac{d\nu}{d\mu}(x)\right)d\nu(x)\mbox{ if }\nu\mbox{ absolutely continuous w.r.t. }\mu\\+\infty\;\mbox{otherwise}\end{cases}\\
\chi^2_2(\nu|\mu)&=\begin{cases}\int_{\R^d} \left(\frac{d\nu}{d\mu}(x)-1\right)^2d\mu(x)
  =\|\frac{d\nu}{d\mu}-1\|^2_{L^2(\mu)}\mbox{ if }\nu\mbox{ absolutely continuous w.r.t. }\mu\\+\infty\;\mbox{otherwise}\end{cases}.
\end{align*}
Next, we precise the inequalities that will be discussed in the paper.
\begin{aadef}The probability measure $\mu$ on $\R^d$ is said to satisfy 
\begin{description}
\item[the Poincaré inequality] ${\cal P}(C)$ with constant $C$ if
$$\forall \varphi:\R^d\to\R\;C^1\;\mbox{with a bounded gradient},\;\int_\R\varphi^2(x)d\mu(x)-\left(\int_\R\varphi(x)d\mu(x)\right)^2\leq C\int_\R|\nabla \varphi(x)|^2d\mu(x)$$
  \item[the transport-chi-square inequality] ${\cal T}_\chi(C)$ with constant $C$ if
$$\forall \nu\mbox{ probability measure on $\R^d$},\;W_2(\mu,\nu)\leq \sqrt{C}\chi_2(\nu|\mu).$$
\item[the log-Sobolev inequality] ${\cal LS}(C)$ with constant $C$ if $\forall \varphi:\R^d\to\R\;C^2$ compactly supported,
$$\int_\R\varphi^2(x)\ln(\varphi^2(x))d\mu(x)-\int_\R\varphi^2(x)d\mu(x)\ln\left(\int_\R\varphi^2(x)d\mu(x)\right)\leq C\int_\R|\nabla \varphi(x)|^2d\mu(x).$$
\item[the transport-entropy inequality] ${\cal T}_H(C)$ with constant $C$ if
$$\forall \nu\mbox{ probability measure on $\R^d$},\;W_2(\mu,\nu)\leq \sqrt{C}H(\nu|\mu).$$
\end{description}
\end{aadef}
According to \cite{ov}, the log-Sobolev inequality is stronger than the transport-entropy inequality which is itself stronger than the Poincaré inequality and more precisely ${\cal LS}(C)\Rightarrow{\cal T}_H(C)\Rightarrow {\cal P}(C/2)$. The transport-entropy inequality is strictly weaker than the log-Sobolev inequality (see \cite{cg,go} for examples of one-dimensional probability measures $\mu$ satisfying the transport-entropy inequality but not the log-Sobolev inequality) and is strictly stronger than the Poincaré inequality (see for example \cite{go} Theorem 1.7).

On the other hand, the inequality $x\ln(x)\leq (x-1)+(x-1)^2$ implies $H(\nu|\mu)\leq \chi_2^2(\nu|\mu)$ and therefore ${\cal T}_H(C)\Rightarrow {\cal T}_\chi(C)$. The transport-entropy inequality implies both the transport-chi-square and Poincaré inequalities. The relation between the two latter is therefore a natural question. It turns out that, by an easy adaptation of the linearization argument in \cite{ov}, the transport-chi-square inequality implies the Poincaré inequality. Moreover, in dimension $d=1$, we are able to prove the converse implication so that both inequalities are equivalent.
Last, we prove tensorization of the transport-chi-square inequality.

{\bf Acknowledgement : } I thank Arnaud Guillin for fruitful discussions and in particular for pointing out the implication ${\mathcal T}_\chi(C)\Rightarrow {\mathcal P}(C)$ and the interest of tensorization to me. 
\section{Main results}
\begin{thm}\label{thequiv}$\forall d\geq 1$, ${\mathcal T}_\chi(C)\Rightarrow {\mathcal P}(C)$. Moreover, when $d=1$, ${\mathcal P}(C)\Rightarrow {\mathcal T}_\chi(32C)$ and the transport-chi-square and Poincaré inequalities are equivalent.
\end{thm}
Before proving Theorem \ref{thequiv}, we state our second main result dedicated to the tensorization property of the transport-chi-square inequality. Its proof is postponed in Section \ref{sectens}.
\begin{thm}\label{thtensor}
   If $\mu_1$ and $\mu_2$ are probability measures on $\R^{d_1}$ and $\R^{d_2}$ respectively satisfying ${\cal T}_\chi(C_1)$ and ${\cal T}_\chi(C_2)$, then the measure $\mu_1\otimes\mu_2$ satisfies ${\cal T}_\chi((C_1+C_2(1+\sqrt{(3d_2+2)d_2}))\wedge (C_2+C_1(1+\sqrt{(3d_1+2)d_1})))$.
\end{thm}
\begin{arem}
  According to Proposition 8.4.1 \cite{abcfgmrs}, if $\mu_1$ and $\mu_2$ respectively satisfy ${\cal T}_H(C_1)$ and ${\cal T}_H(C_2)$, then $\mu_1\otimes\mu_2$ satisfies ${\cal T}_H(C_1\vee C_2)$. The constant that we obtain in the tensorization of the transport-chi-square inequality is larger than $C_1\vee C_2$.
\end{arem}
The proof of the one-dimensional implication ${\cal P}(C)\Rightarrow {\mathcal T}_\chi(32C)$ in Theorem \ref{thequiv} relies on the two next propositions, the proof of which are respectively postponed in Sections \ref{pw2FGf} and \ref{pFGfg}. When $d=1$, we denote by $F(x)=\mu((-\infty,x])$ and  $G(x)=\nu((-\infty,x])$ the cumulative distribution functions of the probability measures $\mu$ and $\nu$. The càg pseudo-inverses of $G$ (resp. $F$) is defined by $G^{-1}:]0,1[\ni u\mapsto\inf\{x\in\R:G(x)\geq u\}$ (resp. $F^{-1}(u)=\inf\{x\in\R:G(x)\geq u\}$) and satisfies
\begin{align}\label{invcdf}
\forall x\in\R,\;\forall u\in(0,1),\;x<G^{-1}(u)\Leftrightarrow G(x)<u.
\end{align}
When $\mu$ (resp. $\nu$)  admits a density w.r.t. the Lebesgue measure, this density is denoted by $f$ (resp. $g$). Moreover, the optimal coupling in \eqref{wass} is given by $\gamma=du\circ(F^{-1},G^{-1})^{-1}$ where $du$ denotes the Lebesgue measure on $(0,1)$ so that $W_q^q(\mu,\nu)=\int_0^1(F^{-1}(u)-G^{-1}(u))^qdu$ (see \cite{rr} p107-109). We take advantage of this optimal coupling to work with the cumulative distribution functions and check the following proposition. In higher dimensions, far less is known on the optimal coupling and this is the main reason why we have not been able to check whether the Poincaré inequality implies the transport-chi-square inequality. 
\begin{aprop}\label{w2FGf}
   If a probability measure $\mu$ on the real line admits a positive probability density $f$, then, for any probability measure $\nu$ on $\R$,
  \begin{align}
W_2^2(\mu,\nu)\leq 4\int_\R\frac{(F-G)^2}{f}(x)dx.\label{FGsurf}
  \end{align}
\end{aprop}
\begin{arem}\begin{itemize}
   \item One deduces that $W^2_1(\mu,\nu)\leq 4\int_\R\frac{(F-G)^2}{f}(x)dx$. Notice that since, by \eqref{invcdf} and Fubini's theorem, \begin{align*}
      W_1(\mu,\nu)&=\int_0^1\int_\R1_{\{F^{-1}(u)\leq x<G^{-1}(u)\}}+1_{\{G^{-1}(u)\leq x<F^{-1}(u)\}}dxdu \\&=\int_\R\int_0^11_{\{G(x)<u\leq F(x)\}}+1_{\{F(x)<u\leq G(x)\}}dudx=\int_\R|F(x)-G(x)|dx,
   \end{align*}
the stronger bound $$W^2_1(\mu,\nu)=\left(\int_\R \frac{|F-G|}{\sqrt{f}}\times\sqrt{f}(x)dx\right)^2\leq \int_\R\frac{(F-G)^2}{f}(x)dx$$ is a consequence of the Cauchy-Schwarz inequality.
\item It is not possible to control $\int_\R\frac{(F-G)^2}{f}(x)dx$ in terms of $W_2^2(\mu,\nu)$. Indeed for $f(x)=\frac{1}{2}e^{-|x|}$ and $d\nu(x)=\frac{1}{2}e^{-|x-m|}dx$, one has $W^2_2(\mu,\nu)=m^2$, $G(x)=\frac{e^{x-m}}{2}1_{\{x\leq m\}}+(1-\frac{e^{m-x}}{2})1_{\{x>m\}}$ and for $m>0$,
$$\int_\R\frac{(F-G)^2}{f}(x)dx\geq \int_m^{+\infty}\frac{(F-G)^2}{f}(x)dx=\frac{e^{-m}}{2}(e^m-1)^2.$$
\end{itemize}
\end{arem}
Next, when the probability measure $\mu$ on the real line admits a positive probability density satisfying a tail assumption known to be equivalent to the Poincaré inequality (see Theorem 6.2.2 \cite{abcfgmrs}), we are able to control the right-hand-side of \eqref{FGsurf} in terms of $\chi_2^2(\nu|\mu)$.
\begin{aprop}\label{FGfg}Let $f(x)$ be a positive probability density on the real line with cumulative distribution function $F(x)=\int_{-\infty}^x f(y)dy$ and median $m$ such that
\begin{align}
   b\stackrel{\rm def}{=}\sup_{x\geq m}\int_{x}^{+\infty}f(y)dy\int_m^x\frac{dy}{f(y)}\vee\sup_{x\leq m}\int^{x}_{-\infty}f(y)dy\int_x^m\frac{dy}{f(y)}<+\infty.\label{defb}
\end{align}  Then for any probability density $g$ on the real line with cumulative distribution function $G(x)=\int_{-\infty}^x g(y)dy$,
\begin{equation}
 \int_\R \frac{(F-G)^2}{f}(x)dx\leq 4b \int_\R\frac{(f-g)^2}{f}(x)dx. \label{majw2FGf}  
\end{equation}
\end{aprop}
\begin{arem}
 \begin{itemize}
   \item The combination of these two propositions implies that any probability measure $\mu$ on the real line admitting a positive density $f$ such that $b<+\infty$ satifies ${\cal T}_\chi(16b)$.
\item Proposition \ref{FGfg} is a generalization of the last assertion in Lemma 2.3 \cite{jm} where $f$ is restricted to the class of probability densities $f_\infty$ solving $f_\infty(x)=-A(F_\infty(x))$ on the real line with
$$A:[0,1]\to\R_-\;C^1,\mbox{ negative on $(0,1)$ and s.t. }A(0)=A(1)=0,\;A'(0)<0,\;A'(1)>0.$$
The constant $b$ associated with any such density is finite by the proof of Lemma 2.1 \cite{jm}. Moreover, in order to investigate the long-time behaviour of the solution $f_t$ of the Fokker-Planck equation
$$\partial_t f_t(x)=\partial_{xx}f_t(x)+\partial_x(A'(F_t(x))f_t(x)),\;(t,x)\in [0,+\infty)\times\R$$
to the density $f_\infty$ such that $\int_\R x f_\infty(x)dx=\int_\R x f_0(x)dx$, \cite{jm} first investigates the exponential convergence to $0$ of $\int_\R\frac{(F_t-F_\infty)^2}{f_\infty}(x)dx$ (Lemma 2.8) before dealing with that of $\int_\R\frac{(f_t-f_\infty)^2}{f_\infty}(x)dx$ (Theorem 2.4).
\item Even when $b<+\infty$, it is not possible to control $\int_\R\frac{(f-g)^2}{f}(x)dx$ in terms of $\int_\R \frac{(F-G)^2}{f}(x)dx$. Indeed let $f(x)=\frac{1}{2}e^{-|x|}$ and 
$$\mbox{ for }n\in\N,\;g_n(x)=\sum_{k\leq n}f(x)1_{[k-1,k)}(|x|)+\sum_{k\geq n}\frac{e^{-\frac{|x|}{2}}}{2}1_{[x_k,k+1)}(|x|)$$
where $x_k=k+1-2\ln\left(1+\frac{e-1}{2}e^{-\frac{k+1}{2}}\right)$ belongs to $(k,k+1)$ and is such that $\int_{x_k}^{k+1}e^{-\frac{x}{2}}dx=\int_k^{k+1}e^{-x}dx$. One has, using $\forall y\geq 0,\;\ln(1+y)\geq \frac{y}{1+y}$ by concavity of the logarithm  and $1+\frac{e-1}{2}e^{-\frac{k+1}{2}}\leq \sqrt{e}$ for the inequality, \begin{align*}
   \int_\R\frac{(f-g_n)^2}{f}(x)dx&=2\int_n^{+\infty}\frac{g_n^2}{f}(x)dx-e^{-n}=2\sum_{k\geq n}\ln\left(1+\frac{e-1}{2}e^{-\frac{k+1}{2}}\right)-e^{-n}\\&\geq \frac{(e-1)}{\sqrt{e}}\sum_{k\geq n}e^{-\frac{k+1}{2}}-e^{-n}=(\sqrt{e}+1)e^{-\frac{n+1}{2}}-e^{-n}.
\end{align*}
On the other hand, since for $k\geq n$ and $x\in[k,k+1]$, $1-\frac{e^{-k}}{2}\leq G_n(x)\leq F(x)=1-\frac{e^{-x}}{2}$, 
\begin{align*}
   \int_\R \frac{(F-G_n)^2}{f}(x)dx\leq \sum_{k\geq n}\int_k^{k+1}\frac{(e^{-k}-e^{-x})^2}{e^{-x}} dx=\frac{e^2-2e-1}{e-1}e^{-n}.
\end{align*}
\end{itemize}
\end{arem}

\begin{adem}[of Theorem \ref{thequiv}]
The implication ${\mathcal T}_\chi(C)\Rightarrow {\mathcal P}(C)$ is obtained by linearization of the transport-chi-square inequality ${\cal T}_\chi(C)$. For $\nu_\varepsilon=(1+\varepsilon \phi)\mu$ with $\phi:\R^d\to\R$ a $C^2$ function compactly supported and such that $\int_{\R^d} \phi(x)d\mu(x)=0$, according to \cite{ov} p394, there is a finite constant $K$ not depending on $\varepsilon$ such that
\begin{align*}
   \int_{\R^d} \phi^2(x)d\mu(x)\leq \sqrt{\int_{\R^d}|\nabla \phi(x)|^2d\mu(x)}\times\frac{W_2(\mu,\nu_\varepsilon)}{\varepsilon}+\frac{KW^2_2(\mu,\nu_\varepsilon)}{\varepsilon}.
\end{align*}
When ${\cal T}_\chi(C)$ holds, then $W_2(\mu,\nu_\varepsilon)\leq\varepsilon\sqrt{C\int_{\R^d}\phi^2(x)d\mu(x)}$ and taking the limit $\varepsilon\to 0$, one deduces that
$$\int_{\R^d} \phi^2(x)d\mu(x)\leq \sqrt{\int_{\R^d}|\nabla \phi(x)|^2d\mu(x)}\times\sqrt{C\int_{\R^d} \phi^2(x)d\mu(x)}.$$
This implies $\int_{\R^d} \phi^2(x)d\mu(x)\leq C\int_{\R^d}|\nabla \phi(x)|^2d\mu(x)$. Let now $\varphi,\phi_n:\R^d\to\R$ be $C^2$ functions compactly supported with $\phi_n$ taking its values in $[0,1]$, equal to $1$ on the ball centered at the origin with radius $n$ and $\nabla \phi_n$ bounded by $1$. Taking the limit $n\to\infty$ in the inequality written with $\phi$ replaced by $\varphi_n=\varphi-\phi_n\frac{\int_{\R^d}\varphi(x)d\mu(x)}{\int_{\R^d}\phi_n(x)d\mu(x)}$, one deduces that the Poincaré inequality ${\cal P}(C)$ holds for $\varphi$. The extension to $C^1$ functions $\varphi$ with a bounded gradient is obtained by density. 

To prove the converse implication, we now suppose that $d=1$, $\mu$ satisfies the Poincaré inequality ${\cal P}(C)$ and that $\chi_2(\nu|\mu)<+\infty$. We set $\mu_n=\rho_{n}\star\mu$ and $\nu_n=\rho_{n}\star \nu$ for $n\geq 1$ where \begin{equation}\label{defrho}
   \rho_n(x)=\sqrt{\frac{n}{2\pi}}e^{-\frac{n x^2}{2}}
\end{equation} denotes the density of the centered Gaussian law with variance $1/n$. For $\varphi$ a $C^1$ function on $\R$ with a bounded derivative such that $0=\int_\R\varphi(x)d\mu_{{n}}(x)=\int_\R \rho_{n}\star\varphi(x)d\mu(x)$, one has 
\begin{align*}
   \int_\R\varphi^2(x)d\mu_n(x)&= \int_\R(\rho_{n}\star\varphi^2)(x)-(\rho_{n}\star\varphi)^2(x)d\mu(x)+\int_\R(\rho_{n}\star\varphi)^2(x)d\mu(x)\\
&\leq \int_\R \frac{1}{n}(\rho_{n}\star(\varphi')^2)(x)d\mu(x)+C\int_\R (\rho_{n}\star\varphi')^2(x)d\mu(x)\\
&\leq \frac{1+nC}{n}\int_\R (\rho_{n}\star(\varphi')^2)(x)d\mu(x)=\frac{1+nC}{n}\int_\R (\varphi')^2(x)d\mu_n(x)
\end{align*}
where we used the Poincaré inequalities for the Gaussian density $\rho_{n}$ (\cite{abcfgmrs} Théorème 1.5.1 p10) applied to $\varphi$ and for $\mu$ applied to $\rho_{n}\star\varphi$ for the second inequality then Jensen's inequality. The probability measure $\mu_n$ admits a positive density w.r.t. the Lebesgue measure and satisfies ${\cal P}(\frac{1+nC}{n})$. According to Theorem 6.2.2 \cite{abcfgmrs}, this property is equivalent to the fact that the constant associated with $\mu_n$ through \eqref{defb} is $b_n\leq 2\frac{1+nC}{n}$. Combining Propositions \ref{w2FGf} and \ref{FGfg}, one deduces that
\begin{equation*}
   W_2^2(\mu_n,\nu_n)\leq 32\frac{1+nC}{n}\chi_2^2(\nu_n|\mu_n).
\end{equation*}
To conclude, let us check that $W_2^2(\mu,\nu)\leq \liminf_{n\to\infty}W_2^2(\mu_n,\nu_n)$ and that $\chi_2^2(\nu_n|\mu_n)\leq \chi_2^2(\nu|\mu)$. First, the probability measures $\mu_n$ with c.d.f. $F_n(x)=\mu_n((-\infty,x])$ (resp $\nu_n$ with c.d.f. $G_n(x)=\nu_n((-\infty,x])$) converge weakly to $\mu$ (resp. $\nu$) which ensures that $du$ a.e. on $(0,1)$, $(F_n^{-1}(u),G_n^{-1}(u))$ tends to $(F^{-1}(u),G^{-1}(u))$ as $n\to\infty$. With Fatou lemma, one deduces that $$W_2^2(\mu,\nu)=\int_0^1(F^{-1}(u)-G^{-1}(u))^2du\leq \liminf_{n\to\infty}\int_0^1(F_n^{-1}(u)-G_n^{-1}(u))^2du=\liminf_{n\to\infty}W_2^2(\mu_n,\nu_n).$$ 
On the other hand, by Jensen's inequality,
\begin{align*}
   \chi_2^2(\nu_n|\mu_n)&=\int_\R\left(\frac{\int_\R(\frac{d\nu}{d\mu}(y)-1)\rho_{n}(x-y)d\mu(y)}{\int_\R\rho_{n}(x-y)d\mu(y)}\right)^2\int_\R\rho_{n}(x-z)d\mu(z)dx\\&\leq\int_{\R}\int_\R\left(\frac{d\nu}{d\mu}(y)-1\right)^2\rho_{n}(x-y)d\mu(y)dx=\chi_2^2(\nu|\mu).
\end{align*}\end{adem}
\begin{arem}\label{cw2n}
Since $$W_2^2(\mu_n,\nu_n)\leq \inf_{\gamma<^{\mu}_{\nu}}\int_{\R^3}((x+z)-(y+z))^2d\gamma(x,y)\rho_{n}(z)dz=W_2^2(\mu,\nu),$$
one has $\lim_{n\to\infty}W_2(\mu_n,\nu_n)=W_2(\mu,\nu)$.\\
Moreover, when $\chi_2^2(\nu|\mu)<+\infty$, then interpreting $\mu_n$ and (resp $\nu_n$) as the distribution at time $\frac{1}{n}$ of a Brownian motion initially distributed according to $\mu$ (resp. $\nu$) and using Theorem 1.7 \cite{fj}, one obtains $\lim_{n\to\infty}\chi_2^2(\nu_n|\mu_n)=\chi_2^2(\nu|\mu)$.
\end{arem}
\section{Proof of Proposition \ref{w2FGf}}\label{pw2FGf} 
To prove the proposition, one first needs to express the Wasserstein distance in terms of the cumulative distribution functions $F$ and $G$ instead of their pseudo-inverses :

\begin{alem}\label{w2fg}
\begin{align}
 W_2^2(\mu,\nu)=\int_{\R^2}\left((F(x\wedge y)-G(x\vee y))^++(G(x\wedge y)-F(x\vee y))^+\right)dydx.\label{w2FG}  
\end{align}
\end{alem}
\begin{adem}[of Lemma \ref{w2fg}]
Let us first suppose that $\mu$ admits a positive continuous density $f$ w.r.t. the Lebesgue measure. Using the change of variables $(v,w)=(F(x),F(y))$ for the third equality then the equivalence $w<F(G^{-1}(u))\Leftrightarrow F^{-1}(w)<G^{-1}(u)\Leftrightarrow G(F^{-1}(w))<u$ deduced from \eqref{invcdf}, one obtains
 \begin{align}
 W_2^2(\mu,\nu)&=\int_0^1(G^{-1}(u)-F^{-1}(u))^2du\notag\\&=2\int_{[0,1]}\int_{\R^2}\left(1_{\{F^{-1}(u)\leq x\leq y<F^{-1}(F(G^{-1}(u)))\}}+1_{\{F^{-1}(F(G^{-1}(u)))\leq x\leq y\leq F^{-1}(u)\}}\right)dxdydu\notag\\
&=2\int_{[0,1]^3}\left(1_{\{u\leq v\leq w<F(G^{-1}(u))\}}+1_{\{F(G^{-1}(u))\leq v\leq w\leq u\}}\right)\frac{dvdw}{f(F^{-1}(v))f(F^{-1}(w))}du\notag\\
&=2\int_{0}^1\int_v^1\int_0^1\left(1_{\{G(F^{-1}(w))<u\leq v\}}+1_{\{w\leq u\leq G(F^{-1}(v))\}}\right)du\frac{dwdv}{f(F^{-1}(w))f(F^{-1}(v))}\notag\\
&=2\int_0^1\int_v^{1}\left((v-G(F^{-1}(w)))^++(G(F^{-1}(v))-w)^+\right)\frac{dwdv}{f(F^{-1}(w))f(F^{-1}(v))}\notag\\
&=2\int_\R\int_x^{+\infty}\left((F(x)-G(y))^++(G(x)-F(y))^+\right)dydx.\label{W2FG}
 \end{align}
By symmetry, one deduces that \eqref{w2FG} holds.

In the general case, one approximates $\mu$ and $\nu$ by the probability measure $\mu_n=\rho_{n}\star \mu$ and $\nu_n=\rho_{n}\star \nu$ (see \eqref{defrho} for the definition of $\rho_{n}$) which admit smooth positive densities w.r.t. the Lebesgue measure. Let $F_n(x)=\mu_n((-\infty,x])$ and $G_n(x)=\nu_n((-\infty,x])$ denote the associated c.d.f.. One has $\lim_{n\to\infty}W_2(\mu_n,\nu_n)=W_2(\mu,\nu)$ according to Remark \ref{cw2n}.
Moreover, by the weak convergence of $\mu_n$ to $\mu$ and $\nu_n$ to $\nu$, $dx$ a.e. on $\R$, $(F_n(x),G_n(x))$ tends to $(F(x),G(x))$. Since, by Jensen's inequality,
$$(F_n(x)-G_n(y))^+=\left(\int(F(x-z)-G(y-z))\rho_{n}(z)dz\right)^+\leq \int(F(x-z)-G(y-z))^+\rho_{n}(z)dz,$$ the right-hand-side of \eqref{W2FG} gets smaller when replacing $(F,G)$ by $(F_n,G_n)$ and tends to the expression with $(F,G)$ as $n\to\infty$ by Fatou lemma. Hence \eqref{W2FG} still holds.\end{adem}
\begin{adem}[of Proposition \ref{w2FGf}] One has
\begin{align}
   \int_x^{+\infty}(F(x)-G(y))^+dy=1_{\{F(x)>G(x)\}}\int_x^{G^{-1}(F(x))}(F(x)-G(y))dy&\leq (F(x)-G(x))^+(G^{-1}(F(x))-x).\label{majoqueuefg}
\end{align}
By Fubini's theorem and a similar argument,
\begin{align*}
  \int_\R\int_x^{+\infty}(G(x)-F(y))^+dydx&=\int_\R\int_{-\infty}^{x}(G(y)-F(x))^+dydx\\
&\leq \int_\R(G(x)-F(x))^+(x-G^{-1}(F(x)))dx
\end{align*}
With \eqref{W2FG} and \eqref{majoqueuefg}, then using Cauchy-Schwarz inequality and the change of variables $u=F(x)$, one deduces that when $\mu$ admits a positive density $f$ w.r.t. the Lebesgue measure, then
\begin{align*}
   W_2^2(\mu,\nu)&\leq 2\int_\R|G(x)-F(x)||x-G^{-1}(F(x))|dx\\&\leq 2\left(\int_\R\frac{(G(x)-F(x))^2}{f(x)}dx\right)^{1/2}\times\left(\int_\R(x-G^{-1}(F(x)))^2f(x)dx\right)^{1/2}\\
&={2}\left(\int_\R\frac{(G(x)-F(x))^2}{f(x)}dx\right)^{1/2}\times\left(\int_0^1(F^{-1}(u)-G^{-1}(u))^2du\right)^{1/2}.
\end{align*}
Recognizing that the second factor in the r.h.s. is equal to $W_2(\mu,\nu)$, one concludes that \eqref{majw2FGf} holds as soon as $W_2(\mu,\nu)<+\infty$.
To prove \eqref{majw2FGf} without assuming finiteness of $W_2(\mu,\nu)$, one defines a sequence $(G_n)_n$ of cumulative distribution functions converging pointwise to $G$ by setting
\begin{equation*}
   G_n(x)=\begin{cases}F(x)\wedge \frac{1}{n}\mbox{ if }x<G^{-1}(\frac{1}{n})\\
G(x)\mbox{ if }x\in[G^{-1}(\frac{1}{n}),G^{-1}(\frac{n-1}{n}))\\
F(x)\vee \frac{n-1}{n}\mbox{ if }x\geq G^{-1}(\frac{n-1}{n})
 \end{cases}
\end{equation*}
For $x<G^{-1}(\frac{1}{n})$, $G(x)<\frac{1}{n}$, $|F(x)-G_n(x)|=(F(x)-\frac{1}{n})^+\leq \min(|F(x)-G(x)|,(F(x)-\frac{1}{n+1})^+)\leq |F(x)-G_{n+1}(x)|$. Similarly, for $x\geq G^{-1}\left(\frac{n-1}{n}\right)$, $G(x)\geq \frac{n-1}{n}$, $|F(x)-G_n(x)|=(\frac{n-1}{n}-F(x))^+\leq \min(|F(x)-G(x)|,(\frac{n}{n+1}-F(x))^+)\leq |F(x)-G_{n+1}(x)|$. As a consequence, for fixed $x\in\R$, the sequence $(|G_n(x)-F(x)|)_{n\in\N}$ is non-decreasing and goes to $|G(x)-F(x)|$ as $n\to\infty$. By monotone convergence, one deduces that $\lim_{n\to+\infty}\int_\R\frac{(G_n-F)^2}{f}(x)dx=\int_\R\frac{(G-F)^2}{f}(x)dx.$ 
Moreover,
\begin{equation*}
   G^{-1}_n(u)=\begin{cases}
F^{-1}(u)\wedge G^{-1}(\frac{1}{n})\mbox{ if }u\leq \frac{1}{n}\\
G^{-1}(u)\mbox{ if }u\in(\frac{1}{n},\frac{n-1}{n}]\\
F^{-1}(u)\vee G^{-1}(\frac{n-1}{n})\mbox{ if }u>\frac{n-1}{n}
\end{cases}.
\end{equation*}
As a consequence, denoting by $\nu_n$ the probability measure with c.d.f. $G_n$, $W_2^2(\mu,\nu_n)=\int_0^1(F^{-1}(u)-G_n^{-1}(u))^2du<+\infty$ and $W_2^2(\mu,\nu)\leq\liminf_{n\to\infty} W_2^2(\mu,\nu_n)$ by Fatou Lemma. One concludes by taking the limit $n\to+\infty$ in  \eqref{majw2FGf} written with $(\nu_n,G_n)$ replacing $(\nu,G)$.
\end{adem}
\section{Proof of Proposition \ref{FGfg}}\label{pFGfg} Let us assume that $b<+\infty$ and $\int_\R\frac{(f-g)^2}{f}(x)dx<+\infty$. By integration by parts, for $n\in\N^*$,
\begin{align}
   \int_{-n}^n \frac{(F-G)^2}{f}(x)dx=\bigg[(F-G)^2(x)\int_m^x\frac{dy}{f(y)}\bigg]_{-n}^{+n}-2\int_{-n}^n (F-G)(f-g)(x)\int_m^x\frac{dy}{f(y)}dx.\label{ip}
\end{align}
For $x$ larger than the median $m$ of the density $f$, by definition of $b$, then by the equality $(F-G)(x)=\int_x^\infty(g-f)(y)dy$ and Cauchy-Schwarz inequality, one has
\begin{align*}
 0\leq (F-G)^2(x)\int_m^x\frac{dy}{f(y)}\leq b\frac{(F-G)^2(x)}{\int_{x}^{+\infty}f(y)dy}=b\frac{\left(\int_x^\infty (f-g)(y)dy\right)^2}{\int_{x}^{+\infty}f(y)dy}\leq b\int_x^\infty \frac{(f-g)^2}{f}(y)dy.
\end{align*}
where the right-hand-side tends to $0$ as $x\to+\infty$ by integrability of $\frac{(f-g)^2}{f}$ on the real line. Similarly, $\lim_{x\to-\infty}(F-G)^2(x)\int_x^m\frac{dy}{f(y)}=0$. Taking the limit $n\to\infty$ in \eqref{ip} and using again the definition of $b$, one deduces that 
\begin{align}
   \int_\R &\frac{(F-G)^2}{f}(x)dx\leq 2b\int_\R |(F-G)(f-g)|(x)\left(\frac{1_{\{x\geq m\}}}{\int_x^\infty f(y)dy}+\frac{1_{\{x<m\}}}{\int^x_{-\infty} f(y)dy}\right)dx.\label{majtrans1}
\end{align}
The product $|(F-G)(f-g)|(x)\times\left(\frac{1_{\{x\geq m\}}}{\int_x^\infty f(y)dy}+\frac{1_{\{x<m\}}}{\int^x_{-\infty} f(y)dy}\right)$ is locally integrable on $\R$ since the first factor is integrable and the second one is locally bounded. Let $a_n<+\infty$ denote the integral of this function on $[-n,n]$.

By Cauchy Schwarz inequality,
\begin{align}
a_n\leq \sqrt{\int_\R\frac{(f-g)^2}{f}(x)dx}\left(\int_{-n}^n f(F-G)^2(x)\left(\frac{1_{\{x\geq m\}}}{\int_x^\infty f(y)dy}+\frac{1_{\{x<m\}}}{\int^x_{-\infty} f(y)dy}\right)^2dx\right)^{1/2}.\label{csan}
\end{align}
Now, setting $\varepsilon_n=\frac{(F-G)^2(n)}{\int_n^\infty f(y)dy}+\frac{(F-G)^2(-n)}{\int^{-n}_{-\infty} f(y)dy}$, we obtain by integration by parts that for $n\geq|m|$, 
\begin{align*}
   &\int_{-n}^n f(F-G)^2(x)\left(\frac{1_{\{x\geq m\}}}{\int_x^\infty f(y)dy}+\frac{1_{\{x<m\}}}{\int^x_{-\infty} f(y)dy}\right)^2dx\\
&=\left[\frac{(F-G)^2(x)}{\int_x^\infty f(y)dy}\right]_m^{n}-2\int_m^{n}\frac{(F-G)(f-g)(x)}{\int_x^\infty f(y)dy}dx-\left[\frac{(F-G)^2(x)}{\int^x_{-\infty} f(y)dy}\right]^m_{-n}+2\int^m_{-n}\frac{(F-G)(f-g)(x)}{\int^x_{-\infty}f(y)dy}dx\\
&=-4(F-G)^2(m)+\varepsilon_n-2\int_{-n}^n(F-G)(f-g)(x)\left(\frac{1_{\{x\geq m\}}}{\int_x^\infty f(y)dy}-\frac{1_{\{x<m\}}}{\int^x_{-\infty} f(y)dy}\right)dx\\
&\leq 2a_n+\varepsilon_n.
\end{align*}
Plugging this estimation in \eqref{csan}, one deduces that
$$\forall n\geq |m|,\;a_n\leq 1_{\{a_n>0\}}\left(2+\frac{\varepsilon_n}{a_n}\right)\int_\R\frac{(f-g)^2}{f}(x)dx.$$
Using that, according to the analysis of the boundary terms in the first integration by parts performed in the proof, $\lim_{n\to+\infty}\varepsilon_n=0$ and that $(a_n)_n$ is non-decreasing, one may take the limit $n\to\infty$ in this inequality to obtain
$$\int_\R |(F-G)(f-g)|(x)\left(\frac{1_{\{x\geq m\}}}{\int_x^\infty f(y)dy}+\frac{1_{\{x<m\}}}{\int^x_{-\infty} f(y)dy}\right)dx\leq 2\int_\R \frac{(f-g)^2}{f}(x)dx.$$
One easily concludes with \eqref{majtrans1}.
\section{Proof of Theorem \ref{thtensor}}\label{sectens}
Let $\nu$ be a probability measure on $\R^{d_1}\times\R^{d_2})$ with respective marginals $\nu_1$ and $\nu_2$ and such that $\chi_2(\nu|\mu_1\otimes\mu_2)<+\infty$, $\rho$ denote the Radon-Nykodym derivative $\frac{d\nu}{d\mu_1\otimes\mu_2}$ and for $x_1\in\R^{d_1}$, $\rho_1(x_1)=\int_{\R^{d_2}}\rho(x_1,x_2)d\mu_2(x_2)$. Notice that $$\chi_2^2(\nu,\mu_1\otimes\mu_2)=\int_{\R^{d_1+d_2}}\left(\rho(x_1,x_2)-1\right)^2d\mu_1(x_1)d\mu_2(x_2).$$

According to the tensorization property of transport costs (see for instance Proposition A.1 \cite{gole}),

\begin{equation}
   W_2^2(\mu_1\otimes\mu_2,\nu)\leq W_2^2(\mu_1,\nu_1)+\int_{\R^{d_1}}1_{\{\rho_1(x_1)>0\}}W_2^2\left(\mu_2,\frac{\rho(x_1,.)}{\rho_1(x_1)}\mu_2\right)d\nu_1(x_1)\label{tensw2}
\end{equation}
By the inequality ${\cal T}_{\chi}(C_1)$ satisfied by $\mu_1$, the equality $\frac{d\nu_1}{d\mu_1}(x_1)=\rho_1(x_1)=\int_{\R^{d_2}}\rho(x_1,x_2)d\mu_2(x_2)$ and Jensen's inequality, one has
\begin{align}
   W_2^2(\mu_1,\nu_1)\leq C_1\chi^2_2(\nu_1|\mu_1)=C_1\int_{\R^{d_1}}(\rho_1(x_1)-1)^2d\mu_1(x_1)\leq C_1\chi_2^2(\nu,\mu_1\otimes\mu_2)\label{w2marg1}.\end{align}
So the first term of the right-hand-side of \eqref{tensw2} is controled by $\chi_2^2(\nu,\mu_1\otimes\mu_2)$.
By the inequality ${\cal T}_{\chi}(C_2)$ satisfied by $\mu_2$, when $\rho_1(x_1)>0$, $W_2^2\left(\mu_2,\frac{\rho(x_1,.)}{\rho_1(x_1)}\mu_2\right)\leq C_2\int_{\R^{d_2}}\left(\frac{\rho(x_1,x_2)}{\rho_1(x_1)}-1\right)^2d\mu_2(x_2)$. Unfortunately, there is no hope to control \begin{align*}
   \int_{\R^{d_1+d_2}}1_{\{\rho_1(x_1)>0\}}\left(\frac{\rho(x_1,x_2)}{\rho_1(x_1)}-1\right)^2&d\nu_1(x_1)d\mu_2(x_2)\\&=\int_{\R^{d_1+d_2}}1_{\{\rho_1(x_1)>0\}}\left(\frac{\rho(x_1,x_2)}{\rho_1(x_1)}-1\right)^2\rho_1(x_1)d\mu_1(x_1)d\mu_2(x_2)
\end{align*}  in terms of $\chi_2^2(\nu,\mu_1\otimes\mu_2)$ because of the possible very small values of $\rho_1(x_1)$. Therefore it is not enough to plug the latter inequality into the right-hand-side of \eqref{tensw2} to conclude that $\mu_1\otimes\mu_2$ satisfies a transport-chi-square inequality. So we are only going to use this inequality for $\rho_1(x_1)\geq \frac{1}{\alpha}$ where $\alpha$ is some constant larger than $1$ to be optimized at the end of the proof. Using Lemma \ref{rhogd} below with $\beta=\alpha$, one obtains
\begin{align}
   \int_{\R^{d_1}}W_2^2\left(\mu_2,\frac{\rho(x_1,.)}{\rho_1(x_1)}\mu_2\right)&1_{\{\rho_1(x_1)\geq \frac{1}{\alpha}\}}d\nu_1(x_1)\notag
  \\& =\alpha C_2\int_{\R^{d_1+d_2}}(\rho(x_1,x_2)-1)^21_{\{\rho_1(x_1)\geq \frac{1}{\alpha}\}}d\mu_1(x_1)d\mu_2(x_2).\label{w2marg2rhogd}
\end{align}
For small positive values of $\rho_1$, we use the estimation of $W_2^2\left(\mu_2,\frac{\rho(x_1,.)}{\rho_1(x_1)}\mu_2\right)$ deduced from the optimal coupling for the total variation distance. 
If $\nu\neq \mu$, let $\varepsilon$ denote a Bernoulli random variable with parameter $p=\int_{\R^{d_2}}\left(\frac{\rho(x_1,x_2)}{\rho_1(x_1)}\wedge 1\right)d\mu_2(x_2)$ and $(X,Y,Z)$ denote an independent $\R^{d_2}\times\R^{d_2}\times\R^{d_2}$-valued random vector with $X$, $Y$ and $Z$ respectively distributed according to $\frac{1}{p}\left(\frac{\rho(x_1,x_2)}{\rho_1(x_1)}\wedge 1\right)d\mu_2(x_2)$, $\frac{1}{1-p}\left(1-\frac{\rho(x_1,x_2)}{\rho_1(x_1)}\right)^+d\mu_2(x_2)$ and $\frac{1}{1-p}\left(\frac{\rho(x_1,x_2)}{\rho_1(x_1)}-1\right)^+d\mu_2(x_2)$. The random variables $\varepsilon X+(1-\varepsilon) Y$ and $\varepsilon X+(1-\varepsilon) Z$ are respectively distributed according to $d\mu_2(x_2)$ and $\frac{\rho(x_1,x_2)}{\rho_1(x_1)}d\mu_2(x_2)$. As a consequence,
\begin{align*}
   W_2^2\left(\mu_2,\frac{\rho(x_1,.)}{\rho_1(x_1)}\mu_2\right)&\leq \E\left((1-\varepsilon)^2(Y-Z)^2\right)=(1-p)\E\left((Y-Z)^2\right)\\&\leq 2(1-p)\left[\E\left(\left(Y-\int_{\R^{d_2}}y_2d\mu_2(y_2)\right)^2\right)+\E\left(\left(Z-\int_{\R^{d_2}}y_2d\mu_2(y_2)\right)^2\right)\right]\\&\leq 2\int_{\R^{d_2}}\left|x_2-\int_{\R^{d_2}}y_2d\mu_2(y_2)\right|^2\left|\frac{\rho(x_1,x_2)}{\rho_1(x_1)}-1\right|d\mu_2(x_2).\end{align*}
One deduces
\begin{align*}
  \int_{\R^{d_1}}&1_{\{0<\rho_1(x_1)<\frac{1}{\alpha}\}}W_2^2\left(\mu_2,\frac{\rho(x_1,.)}{\rho_1(x_1)}\mu_2\right)d\nu_1(x_1)\\&\leq 2\int_{\R^{d_1+d_2}}\left|x_2-\int_{\R^{d_2}}y_2d\mu_2(y_2)\right|^2\left|\rho(x_1,x_2)-\rho_1(x_1)\right|1_{\{\rho_1(x_1)<\frac{1}{\alpha}\}}d\mu_1(x_1)d\mu_2(x_2)\\
&\leq 2\left(\int_{\R^{d_1+d_2}}\left|x_2-\int_{\R^{d_2}}y_2d\mu_2(y_2)\right|^41_{\{\rho_1(x_1)<\frac{1}{\alpha}\}}d\mu_1(x_1)d\mu_2(x_2)\right)^{1/2}\\&\phantom{\leq 2}\times\left(\int_{\R^{d_1+d_2}}(\rho(x_1,x_2)-\rho_1(x_1))^21_{\{\rho_1(x_1)<\frac{1}{\alpha}\}}d\mu_1(x_1)d\mu_2(x_2)\right)^{1/2}\\
&\leq 2C_2\sqrt{(3d_2+2)d_2}\left(\int_{\R^{d_1}}\frac{\alpha^2(\rho_1(x_1)-1)^2}{(\alpha-1)^2}1_{\{\rho_1(x_1)<\frac{1}{\alpha}\}}d\mu_1(x_1)\right)^{1/2}\\&\phantom{2\sqrt{(3d_2+2)d_2}C_2}\times\left(\int_{\R^{d_1+d_2}}[(\rho(x_1,x_2)-1)^2-(\rho_1(x_1)-1)^2]1_{\{\rho_1(x_1)<\frac{1}{\alpha}\}}d\mu_1(x_1)d\mu_2(x_2)\right)^{1/2}\\
&\leq\frac{C_2\alpha \sqrt{(3d_2+2)d_2}}{\alpha-1}\int_{\R^{d_1+d_2}}(\rho(x_1,x_2)-1)^21_{\{\rho_1(x_1)<\frac{1}{\alpha}\}}d\mu_1(x_1)d\mu_2(x_2),\end{align*}
where we used Cauchy Schwarz inequality for the second inequality, then Lemma \ref{mompoinc} below and an explicit computation of the third factor for the third inequality and last the inequality $\sqrt{b}\sqrt{a-b}\leq \frac{a}{2}$ for any $a\geq b\geq 0$. 

Inserting this estimation together with \eqref{w2marg1} and \eqref{w2marg2rhogd} into \eqref{tensw2}, one obtains
\begin{align*}
   W_2^2(\mu_1\otimes\mu_2,\nu)\leq &C_1\chi_2^2(\nu_1,\mu_1)+C_2\alpha\left(1\vee \frac{\sqrt{(3d_2+2)d_2}}{\alpha-1}\right)\chi_2^2(\nu,\mu_1\otimes\mu_2).
\end{align*}
For the optimal choice $\alpha=1+\sqrt{(3d_2+2)d_2}$, one concludes that the measure $\mu_1\otimes\mu_2$ satisfies ${\cal T}_\chi(C_1+C_2(1+\sqrt{(3d_2+2)d_2}))$. Exchanging the roles of $\mu_1$ and $\mu_2$ in the above reasonning, one obtains that $\mu_1\otimes\mu_2$ also satisfies ${\cal T}_\chi(C_2+C_1(1+\sqrt{(3d_1+2)d_1}))$.

\begin{alem}\label{rhogd}
   For $\beta\geq\alpha >0$,
\begin{align*}
   \int_{\R^{d_1+d_2}}\left(\frac{\rho(x_1,x_2)}{\rho_1(x_1)}-1\right)^21_{\{\rho_1(x_1)\geq \frac{1}{\alpha}\}}&d\nu_1(x_1)d\mu_2(x_2)+\beta\int_{\R^{d_1}}(\rho_1(x_1)-1)^21_{\{\rho_1(x_1)\geq \frac{1}{\alpha}\}}d\mu_1(x_1)\\&\leq \beta\int_{\R^{d_1+d_2}}\left(\rho(x_1,x_2)-1\right)^21_{\{\rho_1(x_1)\geq \frac{1}{\alpha}\}}d\mu_1(x_1)d\mu_2(x_2).
\end{align*}
\end{alem}
\begin{adem}
Developping the squares and using the definition of $\rho_1$ and the equality $d\nu_1(x_1)=\rho_1(x_1)d\mu_1(x_1)$, one checks that the difference between the right-hand-side and the first term of the left-hand-side is equal to 
\begin{align*}
   \int_{\R^{d_1}}\left[\left(\beta-\frac{1}{\rho_1(x_1)}\right)\int_{\R^{d_2}}\rho^2(x_1,x_2)d\mu_2(x_2)+(1-2\beta)\rho_1(x_1)+\beta\right]1_{\{\rho_1(x_1)\geq \frac{1}{\alpha}\}}d\mu_1(x_1).
\end{align*}
One easily concludes by remarking that the first integral is retricted to the $x_1\in\R^{d_1}$ such that $\frac{1}{\rho_1(x_1)}\leq \alpha\leq \beta$ and that $\int_{\R^{d_2}}\rho^2(x_1,x_2)d\mu_2(x_2)\geq\left(\int_{\R^{d_2}}\rho(x_1,x_2)d\mu_2(x_2)\right)^2=\rho_1^2(x_1)$.
\end{adem}
\begin{alem}\label{mompoinc}If a probability measure $\mu$ on $\R^d$ satisfies ${\cal T}(C)$, then 
$$\int_{\R^{d}}\left|x-\int_{\R^{d}}yd\mu(y)\right|^2d\mu(x)\leq dC\mbox{ and }\int_{\R^{d}}\left|x-\int_{\R^{d}}yd\mu(y)\right|^4d\mu(x)\leq (3d+2)dC^2.$$
\end{alem}
\begin{adem}
According to Theorem \ref{thequiv}, $\mu$ satisfies ${\cal P}(C)$. By spatial translation, one may assume that $\int_{\R^{d}}yd\mu(y)=0$. Applying the Poincaré inequality ${\cal P}(C)$ to the functions $x=(x_1,\hdots,x_d)\in\R^d\mapsto x_i$, $x\mapsto x_i^2$ and $x\mapsto x_ix_j$ with $1\leq i\neq j\leq d$, yields,
\begin{align*}
   \int_{\R^d}x_i^2d\mu(x)&\leq C\\
\int_{\R^d}x_i^4d\mu(x)&\leq 4C\int_{\R^d}x_i^2d\mu(x)+\left(\int_{\R^d}x_i^2d\mu(x)\right)^2\leq 5C^2\\
\int_{\R^d}(x_ix_j)^2d\mu(x)&\leq C\int_{\R^d}x_i^2+x_j^2d\mu(x)+\left(\int_{\R^d}x_ix_jd\mu(x)\right)^2\leq 2C^2+\int_{\R^d}x_i^2d\mu(x)\int_{\R^d}x_j^2d\mu(x)\leq 3C^2.\end{align*}
One easily concludes by summation of these inequalities.
\end{adem}

\end{document}